\documentclass[12pt]{article}

\usepackage{amsfonts}
\usepackage{amsmath}
\usepackage{amssymb}
\usepackage{color}
\usepackage{amsthm}
\usepackage{hyperref}
\usepackage{makecell}
\usepackage{slashbox,pict2e}

\newtheorem{theorem}{Theorem}
\newtheorem{corollary}[theorem]{Corollary}

\newtheorem{remark}[theorem]{Remark}

\numberwithin{equation}{section}

\title{Poisson and Gaussian approximations of the power divergence family of statistics} \author{Fraser Daly\footnote{Department of Actuarial Mathematics and Statistics and the Maxwell Institute for Mathematical Sciences, Heriot--Watt University, Edinburgh EH14 4AS, UK.  E-mail: f.daly@hw.ac.uk}} \date{\today}

\begin{document}

\maketitle

\noindent{\bf Abstract} 
Consider the family of power divergence statistics based on $n$ trials, each leading to one of $r$ possible outcomes. This includes the log-likelihood ratio and Pearson's statistic as important special cases. It is known that in certain regimes (e.g., when $r$ is of order $n^2$ and the allocation is asymptotically uniform as $n\to\infty$) the power divergence statistic converges in distribution to a linear transformation of a Poisson random variable. We establish explicit error bounds in the Kolmogorov (or uniform) metric to complement this convergence result, which may be applied for any values of $n$, $r$ and the index parameter $\lambda$ for which such a finite-sample bound is meaningful. We further use this Poisson approximation result to derive error bounds in Gaussian approximation of the power divergence statistics.  
\vspace{12pt}


\vspace{12pt}


\section{Introduction and main result}\label{sec:intro}

In the setting of a multinomial goodness-of-fit test, suppose we have $n$ independent trials, each of which is assigned to one of $r=r(n)$ classes. We denote by $\mathbf{p}=(p_1,\ldots,p_r)$ the corresponding classification probabilities, where $p_1+\cdots+p_r=1$, and by $N_1,\ldots,N_r$ the observed frequencies in each class following these $n$ trials, so that $N_1+\cdots+N_r=n$.

Cressie and Read \cite{cressie84} introduced the power divergence family of statistics for testing the hypothesis $H_0: \mathbf{p}=\mathbf{p_0}$ against alternatives of the form $H_1: \mathbf{p}=\mathbf{p_1}$. For a given index parameter $\lambda\in\mathbb{R}$, these statistics are defined by
\[
T_\lambda=n\sum_{j=1}^rp_jg_\lambda\left(\frac{N_j}{np_j}\right)\,,
\]
where 
\begin{equation}\label{eq:gdef}
g_\lambda(x)=\left\{\begin{array}{ll}
2x\log(x)\,, & \text{if }\lambda=0\,,\\
-2\log(x)\,, & \text{if }\lambda=-1\,,\\
\frac{2x(x^\lambda-1)}{\lambda(\lambda+1)}\,, & \text{otherwise}\,.
\end{array}\right.
\end{equation}
As special cases, these statistics include Pearson's $\chi^2$ ($\lambda=1$), the log-likelihood ratio ($\lambda=0$), the Freeman--Tukey statistic ($\lambda=-1/2$) and the modified log-likelihood ratio ($\lambda=-1$). Cressie and Read \cite{cressie84} also suggest $T_{2/3}$ as a good alternative to Pearson's statistic and the log-likelihood ratio. The choice of the index parameter $\lambda$ for applications will depend on the type of departure from $H_0$ that the practitioner wishes to detect. For example, Cressie and Read \cite{cressie84} note that the effect of a single ratio $N_j/p_j$ on the test statistic $T_\lambda$ increases as $|\lambda|$ increases; see their Section 6.2. See also Sections 4.5, 5.5 and 6.7 of \cite{read88} for discussions of appropriate choices of $\lambda$ for testing different hypotheses and in different settings.   

It is well known that in many cases $T_\lambda$ has an asymptotic chi-square distribution as $n\to\infty$. Rates of convergence and error bounds in the chi-square approximation of $T_\lambda$ for $\lambda>-1$ have been obtained to complement and quantify this limit theorem by several authors: see, for example, \cite{assylebekov10,assylebekov11,gaunt22,puchkin23,ulyanov09} and references therein. A Gaussian limit for a suitably scaled version of $T_\lambda$ in some regimes is also well known \cite{cressie84,morris75}. Though we are not aware of any error bounds that accompany this central limit theorem, there are rates of convergence known in particular examples, such as for Pearson's statistic (see Section 4.2 of \cite{dobler23}). Other limiting distributions are also possible. Recently, Rempa\l{}a and Weso\l{}owski \cite{rempala23} have shown that in some regimes, a suitably normalised version of $T_\lambda$ has an asymptotic Poisson distribution as $n\to\infty$; we discuss their result in more detail below. Although this limit theorem is a recent observation, Poisson limits in this setting are not unexpected: Steck \cite{steck57} observed many years ago that in the special case of equally likely classes, Pearson's statistic has an asymptotic Poisson distribution. 

Our main aim in this paper is to complement the limit theorem of \cite{rempala23} for power divergence statistics by providing an explicit error bound in the Poisson approximation of $T_\lambda$ in Theorem \ref{thm:main} below. That is, we give an upper bound on the Kolmogorov (or uniform) distance between $\widetilde{T}_\lambda$, a suitably normalised version of $T_\lambda$ to be defined precisely in \eqref{eq:tilde_def} below, and a suitable Poisson random variable $Z$, given by
\begin{equation}\label{eq:dk}
d_K(\widetilde{T}_\lambda,Z)=\sup_{y\in\mathbb{R}}|\mathbb{P}(\widetilde{T}_\lambda\geq y)-\mathbb{P}(Z\geq y)|\,.
\end{equation}
We will also use this Poisson approximation result to derive error bounds in Gaussian approximation of the power divergence family of statistics, again in the Kolmogorov distance.

This choice of metric is governed by practical considerations in the use of the statistics $T_\lambda$. Theorem \ref{thm:main} can be used to give upper bounds on the probability of type-I and type-II errors in the test of the hypotheses $H_0$ versus $H_1$ described above. For example, to bound the probability of a type-I error, set $\mathbf{p}=\mathbf{p_0}$. One is then interested in bounding the probability that the observed value of $T_\lambda$ falls into the rejection region. This can be done by approximating $T_\lambda$ by a suitably scaled Poisson random variable, for which the corresponding error probability can be calculated or bounded, and bounding the error in this approximation using Theorem \ref{thm:main}. These estimates will be useful in settings where this error estimate is small, and the Poisson approximation is known to hold. Similar comments apply to the estimation of type-II error probabilities. 

In a Poisson approximation setting it is relatively straightforward to directly derive error bounds in the Kolmogorov distance, but this is typically not the case for chi-square and Gaussian approximations. In the setting of chi-square approximation for $T_\lambda$, \cite{gaunt22} was able to directly obtain bounds only in a metric based on smooth test functions. Similarly, the bounds given by \cite{dobler23} in the Gaussian approximation of Pearson's $\chi^2$ are for smooth test functions. Bounds in Kolmogorov distance can be obtained from such results for smooth test functions, but typically such bounds give sub-optimal rates of convergence. Rates of convergence (but not explicit error bounds) in Kolmogorov distance for the chi-square approximation of $T_\lambda$ have also been established directly by several authors: see \cite{assylebekov10,assylebekov11,ulyanov09}, among others. The recent paper \cite{puchkin23} also derives bounds in the Kolmogorov distance, employing some auxiliary randomisation and showing that the upper bounds hold with high probability with respect to this randomisation. 

We will focus here on the case $\lambda>-1$, since this is the largest class of power divergence statistics we can reasonably consider. For $\lambda\leq-1$, if one or more of $N_1,\ldots,N_r$ is zero then we have $T_\lambda=\infty$. Since this will happen with high probability in our regimes of interest (in which a typical application will have $r$ of order $n^2$), bounds on $d_K(\widetilde{T}_\lambda,Z)$ become meaningless in this case.

We use the remainder of this section to discuss the Poisson limit theorem of Rempa\l{}a and Weso\l{}owski \cite{rempala23}, and to state our main Poisson result in Theorem \ref{thm:main} below. Some examples and illustrations of this result are then given in Section \ref{sec:egs}. In Section \ref{sec:normal} we derive our Gaussian approximation result. The proof of Theorem \ref{thm:main} is given in Section \ref{sec:proof}, following which some concluding remarks make up Section \ref{sec:conc}.

Let the random variable $P$ be defined by $\mathbb{P}(P=p_j)=p_j$, for $j=1,\ldots,r$. With this notation, Theorem 2.1 of Rempa\l{}a and Weso\l{}owski \cite{rempala23} establishes a Poisson limit theorem for $T_\lambda$ for any $\lambda\geq0$ under the following three conditions:
\begin{enumerate}
\item[(i)] $n\max_jp_j\to0$ as $n\to\infty$;
\item[(ii)]$\frac{n\text{Var}(P^{-\lambda})}{\mathbb{E}[P^{-\lambda}]^2}\to0$ as $n\to\infty$ when $\lambda>0$, or $n\text{Var}(\log(P))\to0$ as $n\to\infty$ when $\lambda=0$; and
\item[(iii)] there exists a constant $0<\eta<\infty$ such that $\frac{n^2\mathbb{E}[P^{1-k\lambda}]}{\mathbb{E}[P^{-\lambda}]^k}\to\eta$ as $n\to\infty$, for each $k=1,2,3$.
\end{enumerate}
More precisely, under the conditions (i)--(iii) above, Rempa\l{}a and Weso\l{}owski show that the normalised version $\widetilde{T}_\lambda$ of $T_\lambda$, defined by
\begin{equation}\label{eq:tilde_def}
\widetilde{T}_\lambda=\frac{n^\lambda}{g_\lambda(2)\mathbb{E}[P^{-\lambda}]}\left(T_\lambda-n^2\mathbb{E}\left[Pg_\lambda\left([nP]^{-1}\right)\right]\right)\,,
\end{equation}
converges in distribution to a Poisson distribution with mean $\eta/2$ as $n\to\infty$. Our Theorem \ref{thm:main} below gives a quantitative version of this result. Rempa\l{}a and Weso\l{}owski \cite{rempala23} are motivated by \cite{chen10,jankova20,kim20} to consider the doubly asymptotic setting where $r=r(n)\to\infty$ as $n\to\infty$. In considering a rate of convergence based on our Theorem \ref{thm:main} we may also have such a setting in mind, but we emphasise that our result gives an explicit error bound which may be evaluated for any given values of $n$ and $r$.    

\begin{theorem}\label{thm:main}
Let $\lambda>-1$ and $n\geq4$. Assume that $p_j\leq0.13$ and $(n+1)p_j<4$ for each $j=1,\ldots,r$. Let $\pi_j=\binom{n}{2}p_j^2(1-p_j)^{n-2}$ and $\mu=\pi_1+\cdots+\pi_r$, and let $Z$ have a Poisson distribution with mean $\mu$. Then
\begin{multline}\label{eq:main}
d_K(\widetilde{T}_\lambda,Z)\leq45\min\{1,\mu^{-1}\}\left(6n\left[\sum_{j=1}^rp_j\pi_j\right]^2+\sum_{j=1}^r\pi_j^2+\frac{24\mu^2}{n}\right)\\
+18\left(5.55\mu c_\lambda\right)^{0.49/c_\lambda}
+\frac{8n\mu\max_jp_j}{4-(n+1)\max_jp_j}
+\min\left\{8.1d_\lambda,2.15\left(\frac{d_\lambda}{\mu}\right)^{1/3}\right\}\,,
\end{multline}
where
\[
c_\lambda=\max_j\left|\frac{1}{p_j^\lambda\mathbb{E}[P^{-\lambda}]}-1\right|\,,\quad\text{and}\quad
d_\lambda=\left\{\begin{array}{ll}
\frac{n\text{Var}(\log(P))}{4\log(2)^2}\,, & \text{if }\lambda=0\,,\\
\frac{n\text{Var}(P^{-\lambda})}{4(2^\lambda-1)^2\mathbb{E}[P^{-\lambda}]^2}\,, & \text{otherwise}\,.
\end{array}\right.
\]
\end{theorem}

If $c_\lambda=0$, the second term of the upper bound \eqref{eq:main} should be read as zero. The proof of Theorem \ref{thm:main} is given in Section \ref{sec:proof} below, before which we consider some applications of that result in Section \ref{sec:egs} and a Gaussian approximation result which may be derived from it in Section \ref{sec:normal}. We conclude this section with some brief discussion of our upper bound, beginning by noting that while the Poisson limit theorem of \cite{rempala23} applies only for $\lambda\geq0$, our result provides bounds for all $\lambda>-1$. See Section \ref{sec:uniform} below for an example where we obtain a Poisson limit for this full range of values of $\lambda$.

\begin{remark}
The approximating Poisson random variable $Z$ in Theorem \ref{thm:main} has a mean $\mu$ which depends on $n$. In some applications it may be more convenient to replace this with a different Poisson distribution whose mean is equal to $\lim_{n\to\infty}\mu$. Note that under the conditions of Theorem 2.1 of \cite{rempala23}, this limit is equal to $\eta/2$. Letting $Z^\prime$ have a Poisson distribution with mean $\eta/2$, we obtain an upper bound for $d_K(\widetilde{T}_\lambda,Z^\prime)$ by taking the upper bound \eqref{eq:main} of Theorem \ref{thm:main} and adding to it the additional term
\[
45\min\left\{\sqrt{\frac{2}{e}}\left|\sqrt{\mu}-\sqrt{\frac{\eta}{2}}\right|,\left|\mu-\frac{\eta}{2}\right|\right\}\,.
\] 
This is apparent from the proof of Theorem \ref{thm:main} below, by using the triangle inequality at an appropriate point and applying Theorem 1 of \cite{roos03} to bound $d_K(Z,Z^\prime)$, noting that the total variation distance bounded therein serves as an upper bound for the Kolmogorov distance we use here.
\end{remark}

\begin{remark}
The assumptions that $p_j\leq0.13$ and $(n+1)p_j<4$ for each $j$ in Theorem \ref{thm:main} may be relaxed at the cost of a somewhat worse upper bound. However, since we have in mind the setting where $n\max_jp_j\to0$ as $n\to\infty$, we lose very little by making these assumptions here.
\end{remark}

\begin{remark}
It is clear that we need $d_\lambda\to0$ if the final term of \eqref{eq:main} is to go to zero as $n$ grows large, with the rate of convergence then governed by the first term of the minimum therein. However, in particular applications, for a given values of the parameters of the problem, it may be the case that the second term of the minimum is smaller than the first. We therefore leave both terms in the statement of our bound. 
\end{remark}

We emphasise that Theorem \ref{thm:main} applies for any given $\lambda>-1$, $n\geq4$, $r$ and the $p_j$, and can be used to estimate the Poisson approximation error for any set of values of these parameters. Nevertheless, thinking of the limit theorem implied by our upper bound, it is natural to compare it with the conditions (i)--(iii) above for Poisson convergence derived by Rempa\l{}a and Weso\l{}owski \cite{rempala23}. We have that our upper bound \eqref{eq:main} converges to zero as $n\to\infty$ if conditions (i), (ii), (iv) and (v) hold, where these latter two conditions are given as follows:
\begin{enumerate}
\item[(iv)] there exists a constant $0<\eta<\infty$ such that $n^2\mathbb{E}[P]\to\eta$ as $n\to\infty$; and
\item[(v)] $c_\lambda\to0$ as $n\to\infty$.
\end{enumerate} 
To see this, we note that convergence of the final term of our upper bound follows immediately from (ii). Noting that $\mu\leq\frac{n^2}{2}\sum_{j=1}^rp_j^2=\frac{n^2}{2}\mathbb{E}[P]\to\eta/2$ by (iv), convergence of the third term to zero is a consequence of (i). Condition (v) gives convergence of the second term. For the first term of the upper bound, (iv) has already given a finite limit for $\mu$, so that $\mu^2/n\to0$.  Similarly, $\sum_{j=1}^r\pi_j^2\leq\frac{n^4}{4}(\max_jp_j)^2\mathbb{E}[P]\to0$ by (i) and (iv). Finally, 
\[
n\left[\sum_{j=1}^rp_j\pi_j\right]^2\leq n\left[\sum_{j=1}^r\frac{n^2}{2}p_j^3\right]^2\leq\frac{n^5}{4}(\max_jp_j)^2\mathbb{E}[P]^2\to0\,,
\]
again using (i) and (iv).

Some comparison of the conditions (i)--(iii) used in the Poisson limit theorem of \cite{rempala23} and our Theorem \ref{thm:main} will be made in the examples considered in Section \ref{sec:egs}. We remark here that we have $c_\lambda=0$ when either $\lambda=0$ or the allocation is uniform (i.e., $p_j=1/r$ for each $j=1,\ldots,r$). So (v) suggests that the allocation should be asymptotically uniform. In either of these cases where $c_\lambda=0$, we have that (iii) and (iv) are equivalent, since (iii) reduces to the statement that $n^2\mathbb{E}[P]\to\eta$. The proof of Theorem \ref{thm:main} makes it clear that our upper bound is geared towards the natural setting where the allocation is asymptotically uniform and $c_\lambda\to0$, in which case it appears that conditions (iii) and (iv) play similar roles.

\section{Examples and illustrations}\label{sec:egs}

In this section we illustrate the bound \eqref{eq:main} of Theorem \ref{thm:main} with some examples and applications. As noted above, if either $\lambda=0$ or the allocation is uniform then $c_\lambda=0$, and so the second term of the upper bound \eqref{eq:main} vanishes. We will detail both these special cases here, beginning in Section \ref{sec:llr} with the log-likelihood ratio, $\lambda=0$, as an initial illustration. For the case of uniform allocations, we will first study the more general example of the discrete power distribution in Section \ref{sec:dpd}, which we will then specialise to the case of uniform allocation in Section \ref{sec:uniform}.

\subsection{The log-likelihood ratio}\label{sec:llr}

In the case $\lambda=0$ the upper bound \eqref{eq:main} simplifies, giving the following Poisson approximation bound for the log-likelihood ratio $T_0=2\sum_{j=1}^rN_j\log\left(\frac{N_j}{np_j}\right)$.
\begin{corollary}\label{cor:llr}
Let $n\geq4$, and assume that $p_j\leq0.13$ and $(n+1)p_j<4$ for each $j=1,\ldots,r$. Let $\pi_j$ and $\mu$ be as in Theorem \ref{thm:main}, and let $Z$ have a Poisson distribution with mean $\mu$. Then
\begin{multline*}
d_K(\widetilde{T}_0,Z)\leq45\min\{1,\mu^{-1}\}\left(6n\left[\sum_{j=1}^rp_j\pi_j\right]^2+\sum_{j=1}^r\pi_j^2+\frac{24\mu^2}{n}\right)\\
+\frac{8n\mu\max_jp_j}{4-(n+1)\max_jp_j}
+\min\left\{8.1d_0,2.15\left(\frac{d_0}{\mu}\right)^{1/3}\right\}\,,
\end{multline*}
where $d_0=\frac{n\text{Var}(\log(P))}{4\log(2)^2}$.
\end{corollary}
As noted above, since $c_0=0$, the condition (iii) from the limit theorem of Rempa\l{}a and Weso\l{}owski \cite{rempala23} is equivalent to (iv) above in this example.

\subsection{The discrete power distribution}\label{sec:dpd}

We now consider the example in which $p_j=\frac{j^{-a}}{z_r(a)}$ for $j=1,\ldots,r$ and for some parameter $a\in[0,1]$, where the normalising constant is the generalised harmonic number defined by $z_r(a)=\sum_{k=1}^rk^{-a}$. This is the discrete power distribution with parameters $a$ and $r$. This example was also considered in Example 3.1 of \cite{rempala23}, where it is shown that the restriction $a<1/2$ is necessary to ensure that $n\max_jp_j\to0$ as $n\to\infty$ in the case where $r$ is of order $n^2$. In doing so they used the asymptotics $z_r(a)\approx\frac{r^{1-a}}{1-a}$ for $a\not=1$, which we can also employ in understanding the behaviour of the terms that appear in our bound in this setting. For future use we also note that $z_r(1)=O(\log r)$.

In the special case $a=0$ we have the uniform distribution $p_j=1/r$ for each $j$; we will look in more detail at this special case in Section \ref{sec:uniform} below, but begin here with the more general example. For concreteness we will focus on the case $\lambda>0$ throughout this section.  

In this setting we have that $\max_jp_j=p_1=z_r(a)^{-1}$,
\[
\pi_j=\binom{n}{2}\frac{j^{-2a}}{z_r(a)^2}\left(1-\frac{j^{-a}}{z_r(a)}\right)^{n-2}\leq\frac{n^2j^{-2a}}{2z_r(a)^2}\,,\quad\text{and}\quad
\mu=\sum_{j=1}^r\pi_j\leq\sum_{j=1}^r\frac{n^2j^{-2a}}{2z_r(a)^2}=\frac{n^2z_r(2a)}{2z_r(a)^2}\,.
\]
We further note that, for $x\in\mathbb{R}$,
\[
\mathbb{E}[P^x]=\sum_{j=1}^rp_j^{1+x}=\frac{z_r(a(1+x))}{z_r(a)^{1+x}}\,,
\]
and hence
\[
d_\lambda=\frac{n}{4(2^\lambda-1)^2}\left(\frac{\mathbb{E}[P^{-2\lambda}]}{\mathbb{E}[P^{-\lambda}]^2}-1\right)=\frac{n}{4(2^\lambda-1)^2}\left(\frac{z_r(a(1-2\lambda))z_r(a)}{z_r(a(1-\lambda))^2}-1\right)\,.
\]
Finally, we have that
\begin{equation}\label{eq:dpd_c}
c_\lambda=\max_j\left|\frac{z_r(a)j^{a\lambda}}{z_r(a(1-\lambda))}-1\right|=\frac{z_r(a)r^{a\lambda}}{z_r(a(1-\lambda))}-1\,.
\end{equation}
Combining all these ingredients, we have that \eqref{eq:main} gives us
\begin{multline}\label{eq:dpd}
d_K(\widetilde{T}_\lambda,Z)\leq
\frac{45n^3}{z_r(a)^4}\left(\frac{3n^2z_r(3a)^2}{2z_r(a)^2}+\frac{nz_r(4a)}{4}+6z_r(2a)^2\right)
+18\left(\frac{2.78n^2c_\lambda z_r(2a)}{z_r(a)^2}\right)^{0.49/c_\lambda}\\
+\frac{4n^3z_r(2a)}{z_r(a)^2(4z_r(a)-(n+1))}
+\frac{2.025n}{(2^\lambda-1)^2}\left(\frac{z_r(a(1-2\lambda))z_r(a)}{z_r(a(1-\lambda))^2}-1\right)\,,
\end{multline}
where $c_\lambda$ is given by \eqref{eq:dpd_c} and for simplicity we have taken the first term in each of the minima in \eqref{eq:main}. As noted above, we do not lose much by doing this. 

For illustration, consider the case where $r$ is of order $n^2$, as in Example 3.1 of \cite{rempala23}. Of the four terms in the upper bound \eqref{eq:dpd}, the second and fourth rely on $a$ being close to zero in order to be small, whereas the first and third will go to zero as $n\to\infty$ for any fixed $a<1/2$. These terms are of order $O(n^{-1})$ and $O(n^{2a-1})$, respectively, up to logarithmic terms. The bound \eqref{eq:dpd} continues to hold for fixed $a\geq1/2$, but these terms will no longer go to zero as $n$ grows and cannot be expected to be small. In particular, the third term will grow logarithmically in $n$ at $a=1/2$ and as a power of $n$ for $a>1/2$.  

\subsection{Uniform allocations}\label{sec:uniform}

We now specialise the example of Section \ref{sec:dpd} to the case where the allocation is uniform, i.e., $a=0$ and so $p_j=1/r$ for each $j=1,\ldots,r$. Noting that $z_r(0)=r$ and that $c_\lambda=d_\lambda=0$ in this case, we obtain the following corollary to Theorem \ref{thm:main}:
\begin{corollary}\label{cor:uniform}
Let $n\geq4$, $r\geq8$ and $(n+1)/r<4$. Let $T_\lambda$ be the power divergence statistic in the uniform case $p_j=1/r$ for $j=1,\ldots,r$, and $Z$ have a Poisson distribution with mean $\frac{1}{r}\binom{n}{2}(1-1/r)^{n-2}$. Then, for all $\lambda>-1$,
\begin{equation}\label{eq:uniform}
d_K(\widetilde{T}_\lambda,Z)\leq\frac{45n^3}{r^4}\left(\frac{3n^2}{2}+\frac{nr}{4}+6r^2\right)
+\frac{4n^3}{r(4r-(n+1))}\,.
\end{equation}
\end{corollary}
  
It is notable that the upper bound of Corollary \ref{cor:uniform} does not depend on $\lambda$, and holds for all $\lambda>-1$. It follows immediately from the upper bound \eqref{eq:dpd} when $\lambda>0$. Corollary \ref{cor:llr} shows that \eqref{eq:uniform} also holds for $\lambda=0$. In the case $\lambda\in(-1,0)$, we still have that $c_\lambda=d_\lambda=0$, and \eqref{eq:uniform} follows directly from Theorem \ref{thm:main}.

We also note that this converges to zero, with rate $O(1/n)$, when $r$ is of order $n^2$. This is not surprising in light of the conditions (i)--(iii) above for Poisson convergence due to \cite{rempala23}. In this setting of uniform allocations, condition (i) asserts that we must have $n/r\to0$, condition (ii) is trivially satisfied since $d_\lambda=0$, and condition (iii) asserts that that $n^2/r$ has a finite and non-zero limit. We conjecture that $O(1/n)$ is the correct rate of convergence here, noting that the convergence of the underlying multinomial statistics to Poisson, which yields the first term of the upper bound \eqref{eq:uniform}, is known to be of this order in total variation distance in this case; see the discussion on page 121 of \cite{barbour92}. Whether our Theorem \ref{thm:main} provides optimal or near-optimal rates of convergence in other settings is an open problem.  

As well as a rate of convergence, Corollary \ref{cor:uniform} provides an explicit error bound, which we illustrate in Table \ref{tab:uniform} for various choices of $n$ and $r$. As expected from the conditions for Poisson convergence here, this bound is useful only when $r$ is very large compared to $n$. In other regimes a chi-square or Gaussian approximation may be more appropriate than our Poisson approximation. Nevertheless, in such a very sparse regime with very large $r$, the results of Table \ref{tab:uniform} show that our error bound can be small and a Poisson approximation of practical use. We also expect that, while the rate of convergence we obtain is reasonable, there is room for improvement in the leading constants in our error bound, which would improve its usefulness.     

\begin{table}
\begin{center}
\begin{tabular}{r|r|r|r|r|r|r}
\hline
\backslashbox{$r$}{$n$~~}& \multicolumn{1}{c|}{\bf{5}} & \multicolumn{1}{c|}{\bf{10}} & \multicolumn{1}{c|}{\bf{20}} & \multicolumn{1}{c|}{\bf{30}} & \multicolumn{1}{c|}{\bf{40}} & \multicolumn{1}{c}{\bf{50}} \\
\hline
\bf{300}   & 0.3767 & \multicolumn{1}{c|}{---} & \multicolumn{1}{c|}{---} & \multicolumn{1}{c|}{---} & \multicolumn{1}{c|}{---} & \multicolumn{1}{c}{---} \\
\bf{500}   & 0.1356 & \multicolumn{1}{c|}{---} & \multicolumn{1}{c|}{---} & \multicolumn{1}{c|}{---} & \multicolumn{1}{c|}{---} & \multicolumn{1}{c}{---} \\
\bf{700}   & 0.0692 & 0.5534 & \multicolumn{1}{c|}{---} & \multicolumn{1}{c|}{---} & \multicolumn{1}{c|}{---} & \multicolumn{1}{c}{---} \\
\bf{1000}  & 0.0339 & 0.2711 & \multicolumn{1}{c|}{---} & \multicolumn{1}{c|}{---} & \multicolumn{1}{c|}{---} & \multicolumn{1}{c}{---} \\
\bf{3000}  & 0.0038 & 0.0301 & 0.2410 & 0.8134 & \multicolumn{1}{c|}{---} & \multicolumn{1}{c}{---} \\
\bf{5000}  & 0.0014 & 0.0108 & 0.0867 & 0.2928 & 0.6940 & \multicolumn{1}{c}{---} \\
\bf{7000}  & 0.0007 & 0.0055 & 0.0443 & 0.1494 & 0.3540 & 0.6915 \\
\bf{10000} & 0.0003 & 0.0027 & 0.0217 & 0.0732 & 0.1735 & 0.3388 \\
\hline
\end{tabular}
\caption{The upper bound \eqref{eq:uniform} to 4 d.p., for various values of $n$ and $r$. The symbol `---' indicates that the bound is greater than 1, and so not informative.}\label{tab:uniform}
\end{center}
\end{table}

\section{A Gaussian approximation theorem}\label{sec:normal}

Our Poisson approximation result may be used to derive a corresponding error bound in the Gaussian approximation of $\widetilde{T}_\lambda$, as we show here. Letting $\mu$ and $Z$ be as in Theorem \ref{thm:main}, and $N\sim\text{N}(0,1)$ have a standard Gaussian distribution, we use the triangle inequality to write
\[
d_K\left(\frac{\widetilde{T}_\lambda-\mu}{\sqrt{\mu}},N\right)
\leq d_K\left(\frac{\widetilde{T}_\lambda-\mu}{\sqrt{\mu}},\frac{Z-\mu}{\sqrt{\mu}}\right)
+d_K\left(\frac{Z-\mu}{\sqrt{\mu}},N\right)
\leq d_K(\widetilde{T}_\lambda,Z)+\frac{0.4748}{\sqrt{\mu}}\,,
\]
where the final inequality follows from a standard application of the Berry--Esseen inequality, with constant 0.4748 due to \cite{shevtsova11}. This gives the following Gaussian approximation result.

\begin{theorem}\label{thm:normal}
Let $\lambda>-1$ and $n\geq4$. Assume that $p_j\leq0.13$ and $(n+1)p_j<4$ for each $j=1,\ldots,r$. Let $\pi_j$, $\mu$, $c_\lambda$ and $d_\lambda$ be as in Theorem \ref{thm:main}, and let $N\sim\text{N}(0,1)$. Then
\begin{multline*}
d_K\left(\frac{\widetilde{T}_\lambda-\mu}{\sqrt{\mu}},N\right)\leq45\min\{1,\mu^{-1}\}\left(6n\left[\sum_{j=1}^rp_j\pi_j\right]^2+\sum_{j=1}^r\pi_j^2+\frac{24\mu^2}{n}\right)\\
+18\left(5.55\mu c_\lambda\right)^{0.49/c_\lambda}
+\frac{8n\mu\max_jp_j}{4-(n+1)\max_jp_j}
+\min\left\{8.1d_\lambda,2.15\left(\frac{d_\lambda}{\mu}\right)^{1/3}\right\}
+\frac{0.4748}{\sqrt{\mu}}\,.
\end{multline*}
\end{theorem}

With the addition of this final term in our upper bound, we now work in a regime where $\mu$ is large (so that this final term is small), but not so large that the other terms in our upper bound become significant.  We are not aware of other explicit Gaussian error bounds for the family of power divergence statistics with which we can compare our Theorem \ref{thm:normal}, but we will use the remainder of this section to compare our bound with conditions under which convergence to a limiting Gaussian distribution is known to hold. For concreteness, we will do this in the particular example of uniform allocations, as considered in Section \ref{sec:uniform} above. Given that our bound is derived via a Poisson approximation result, it is perhaps not surprising that we will see that conditions under which our upper bound goes to zero are sub-optimal. Nevertheless, we stress again that we are not aware of other bounds that cover the family of power divergence statistics, and that our result gives explicit error bounds in a metric of practical importance.

In the uniform allocation setting with $p_j=1/r$, we have $\pi_j=\binom{n}{2}r^{-2}(1-r^{-1})^{n-2}$ and $c_\lambda=d_\lambda=0$. The upper bound of Theorem \ref{thm:normal} thus gives
\begin{multline}\label{eq:uni_normal}
d_K\left(\frac{\widetilde{T}_\lambda-\mu}{\sqrt{\mu}},N\right)
\leq\frac{90n^3}{(n-1)^2r^3}\left(1-\frac{1}{r}\right)^{2-n}\left(\frac{3n^2}{2}+\frac{nr}{4}+6r^2\right)\\
+\frac{4n^3}{r(4r-(n+1))}+\frac{0.4748}{n-1}\sqrt{2r\left(1-\frac{1}{r}\right)^{2-n}}\,,
\end{multline}  
where, for simplicity of presentation, we write $\mu\leq\frac{n^2}{2r}$ and 
\[
\min\{1,\mu^{-1}\}\leq\mu^{-1}\leq\frac{2r}{(n-1)^2}\left(1-\frac{1}{r}\right)^{2-n}\,.
\]
Recall that in Corollary \ref{cor:uniform} we wrote $\min\{1,\mu^{-1}\}\leq1$, since there we were not in a regime where we expect $\mu$ to be large, hence the differences in the first terms of these two upper bounds for this example.

Considering $r$ of order $n^\alpha$, the upper bound of \eqref{eq:uni_normal} goes to zero as $n\to\infty$ if $\alpha\in(3/2,2)$. This can be compared with the rate of convergence in the Wasserstein distance (based on smooth test functions, unlike the Kolmogorov distance we use here) given by D\"obler \cite{dobler23} for the special case of Pearson's statistic. In his Example 4.6, D\"obler shows that we expect a Gaussian limit for a suitably normalised Pearson's statistic if $\alpha<2$. While our bound is clearly significantly inferior in this respect, it applies to a wider family of statistics and gives an explicit error bound, not only a rate of convergence. 

In the more general setting of the family of power divergence statistics with $\lambda>-1$, Cressie and Read \cite{cressie84} show in their Corollary 2.2 that in the uniform allocation setting with $n/r$ converging to a finite and non-zero limit as $n\to\infty$ a suitably normalised version of $T_\lambda$ has a Gaussian limiting distribution; this fact is not reflected in our upper bound \eqref{eq:uni_normal}. We leave the problem of giving improved error bounds in Kolmogorov distance for the family of power divergence statistics open for future work. 

\section{Proof of Theorem \ref{thm:main}}\label{sec:proof}

We use this section to prove our  main result, the upper bound \eqref{eq:main} of Theorem \ref{thm:main}. Our starting point is the following representation of $T_\lambda$, given in Proposition 2.4 of \cite{rempala23}:
\begin{equation}\label{eq:trep}
T_\lambda=\frac{1}{n^\lambda}\sum_{j=1}^rp_j^{-\lambda}g_\lambda(N_j)+n\sum_{k=1}^np_{X_k}g_\lambda\left([np_{X_k}]^{-1}\right)\,,
\end{equation}
where $g_\lambda$ is given by \eqref{eq:gdef}, and we define the random variable $X$ by $\mathbb{P}(X=j)=p_j$ and let $X,X_1,X_2,\ldots,X_n$ be an IID sequence, so that we can write $N_j=\sum_{k=1}^nI(X_k=j)$. The representation \eqref{eq:trep} follows from the observation that, for any $a,b\geq0$, we have $g_\lambda(ab)=b^{\lambda+1}g_\lambda(a)+ag_\lambda(b)$. Proposition 2.4 of \cite{rempala23} is stated only for $\lambda\geq0$, but the same argument and result hold for any $\lambda>-1$. Combining \eqref{eq:trep} with \eqref{eq:tilde_def}, we write $\widetilde{T}_\lambda=W_\lambda+R_\lambda$, where
\begin{align*}
W_\lambda&=\frac{1}{g_\lambda(2)\mathbb{E}[P^{-\lambda}]}\sum_{j=1}^rp_j^{-\lambda}g_\lambda(N_j)\,,\quad\text{and}\\
R_\lambda&=\frac{n^{\lambda+1}}{g_\lambda(2)\mathbb{E}[P^{-\lambda}]}\left(\sum_{k=1}^np_{X_k}g_\lambda\left([np_{X_k}]^{-1}\right)-n\mathbb{E}\left[Pg_\lambda\left([nP]^{-1}\right)\right]\right)\,.
\end{align*}

Now, for some $\beta_0>0$ to be determined later, we define the event $E_0=\{|R_\lambda|\leq\beta_0\}$ with $\delta_0=\delta_0(\beta_0)=\mathbb{P}(E_0^c)$. Chebyshev's inequality gives
\begin{align}\label{eq:r1_1}
\nonumber\delta_0&=\mathbb{P}\left(\left|\sum_{k=1}^np_{X_k}g_\lambda\left([np_{X_k}]^{-1}\right)-n\mathbb{E}\left[Pg_\lambda\left([nP]^{-1}\right)\right]\right|>\beta_0|g_\lambda(2)|\mathbb{E}[P^{-\lambda}]n^{-(\lambda+1)}\right)\\
&\leq\frac{\text{Var}\left(Pg_\lambda\left([nP]^{-1}\right)\right)n^{2\lambda+3}}{\beta_0^2g_\lambda(2)^2\mathbb{E}[P^{-\lambda}]^2}
=\frac{d_\lambda}{\beta_0^2}\,,
\end{align}
where the final equality follows from the definition of $g_\lambda$, checking the cases $\lambda=0$ and $\lambda\not=0$ separately.

Using the definition \eqref{eq:dk}, we have $d_{K}(\widetilde{T}_\lambda,Z)=\sup_{y\in\mathbb{R}}|\mathbb{P}(\widetilde{T}_\lambda\geq y-\beta_0)-\mathbb{P}(Z\geq y-\beta_0)|$, where the triangle inequality gives
\begin{align}\label{eq:r1_2}
\nonumber|\mathbb{P}(\widetilde{T}_\lambda\geq y-\beta_0)-\mathbb{P}(Z\geq y-\beta_0)|
&\leq|\mathbb{P}(W_\lambda\geq y-\beta_0)-\mathbb{P}(Z\geq y-\beta_0)|+\mathbb{P}(\widetilde{T}_\lambda\geq y-\beta, E_0^c)\\
\nonumber&\qquad+|\mathbb{P}(\widetilde{T}_\lambda\geq y-\beta_0,E_0)-\mathbb{P}(W_\lambda\geq y-\beta_0)|\\
&\leq d_K(W_\lambda,Z)+\delta_0+|\mathbb{P}(\widetilde{T}_\lambda\geq y-\beta_0,E_0)-\mathbb{P}(W_\lambda\geq y-\beta_0)|\,.
\end{align}
For this final term, we write
\begin{align}\label{eq:r1_3}
\nonumber|\mathbb{P}(\widetilde{T}_\lambda\geq y-\beta_0,E_0)&-\mathbb{P}(W_\lambda\geq y-\beta_0)|\\
\nonumber&\leq
|\mathbb{P}(\widetilde{T}_\lambda\geq y-\beta_0,E_0)-\mathbb{P}(W_\lambda\geq y-\beta_0,E_0)|+\mathbb{P}(W_\lambda\geq y-\beta_0,E_0^c)\\
\nonumber&\leq|\mathbb{P}(W_\lambda\geq y-R_\lambda-\beta_0,|R_\lambda|\leq\beta_0)-\mathbb{P}(W_\lambda\geq y-\beta_0,|R_\lambda|\leq\beta_0)|+\delta_0\\
\nonumber&=\mathbb{P}(y-\beta_0-R_\lambda\leq W_\lambda<y-\beta_0,0\leq R_\lambda\leq\beta_0)\\
\nonumber&\qquad+\mathbb{P}(y-\beta_0\leq W_\lambda<y-\beta_0-R_\lambda,-\beta_0\leq R_\lambda<0)+\delta_0\\
\nonumber&\leq\mathbb{P}(y-2\beta_0\leq W_\lambda<y)+\delta_0\\
&\leq2d_{K}(W_\lambda,Z)+\lfloor2\beta_0\rfloor\sup_x\mathbb{P}(Z=x)+\delta_0\,,
\end{align}
where $\lfloor\cdot\rfloor$ is the floor function. 

Combining \eqref{eq:r1_1}--\eqref{eq:r1_3}, we thus have
\begin{equation}\label{eq:r1_4}
d_{K}(\widetilde{T}_\lambda,Z)\leq3d_{K}(W_\lambda,Z)+\frac{2d_\lambda}{\beta_0^2}+\lfloor2\beta_0\rfloor\sup_x\mathbb{P}(Z=x)\,,
\end{equation}
for any $\beta_0>0$. We bound the final two terms on the right-hand side of \eqref{eq:r1_4} in two different ways. Firstly, choosing $\beta_0\approx0.497$, so that $1/\beta_0^2=4.05$ and $\lfloor2\beta_0\rfloor=0$, we have $d_{K}(\widetilde{T}_\lambda,Z)\leq3d_{K}(W_\lambda,Z)+8.1d_\lambda$. Secondly, we may use Proposition A.2.7 of \cite{barbour92} to get that $\mathbb{P}(Z=x)\leq(2e\mu)^{-1/2}$, so that
\[
d_{K}(\widetilde{T}_\lambda,Z)\leq3d_{K}(W_\lambda,Z)+2\left(\frac{d_\lambda}{\beta_0^2}+\frac{\beta_0}{\sqrt{2e\mu}}\right)\,.
\]
We now choose $\beta_0=(2d_\lambda\sqrt{2e\mu})^{1/3}$ to obtain $d_{K}(\widetilde{T}_\lambda,Z)\leq3d_{K}(W_\lambda,Z)+2.15\left(\frac{d_\lambda}{\mu}\right)^{1/3}$, and hence
\begin{equation}\label{eq:r1_final}
d_{K}(\widetilde{T}_\lambda,Z)\leq3d_{K}(W_\lambda,Z)+\min\left\{8.1d_\lambda,2.15\left(\frac{d_\lambda}{\mu}\right)^{1/3}\right\}\,.
\end{equation}

Now, noting that $g_\lambda(0)=g_\lambda(1)=0$, we may write $W_\lambda=V_\lambda+R^\prime_\lambda$, where
\[
V_\lambda=\frac{1}{\mathbb{E}[P^{-\lambda}]}\sum_{j\,:\,N_j=2}p_j^{-\lambda}\,,\quad\text{and}\quad R^\prime_\lambda=\frac{1}{g_\lambda(2)\mathbb{E}[P^{-\lambda}]}\sum_{j\,:\,N_j\geq3}p_j^{-\lambda}g_\lambda(N_j)\,.
\]
We now argue similarly to above: for some $\beta_1>0$ to be determined later, we define the event $E_1=\{R_\lambda^\prime\leq\beta_1\}$ with $\delta_1=\delta_1(\beta_1)=\mathbb{P}(E_1^c)$, noting that $R_\lambda^\prime\geq0$ almost surely for any $\lambda>-1$. We have
\begin{align}\label{eq:r2_1}
\nonumber\delta_1&=\mathbb{P}\left(\frac{1}{g_\lambda(2)\mathbb{E}[P^{-\lambda}]}\sum_{j\,:\,N_j\geq3}p_j^{-\lambda}g_\lambda(N_j)>\beta_1\right)
\leq\mathbb{P}\left(\exists j\in\{1,\ldots,r\}\text{ such that }N_j\geq3\right)\\
&\leq\sum_{j=1}^r\mathbb{P}(N_j\geq3)
\leq4\binom{n}{3}\sum_{j=1}^r\frac{p_j^3(1-p_j)^{n-2}}{4-(n+1)p_j}
\leq\frac{4n\mu\max_jp_j}{3[4-(n+1)\max_jp_j]}\,,
\end{align}
where the penultimate inequality uses Proposition A.2.5(ii) of \cite{barbour92} and the assumption that $(n+1)p_j<4$ for each $j$, since $N_j\sim\text{Bin}(n,p_j)$ has a binomial distribution.

We now write $d_K(W_\lambda,Z)=\sup_{y\in\mathbb{R}}|\mathbb{P}(W_\lambda\geq y-\beta_1)-\mathbb{P}(Z\geq y-\beta_1)|$, where the triangle inequality gives
\begin{align}\label{eq:r2_2}
\nonumber|\mathbb{P}(W_\lambda\geq y-\beta_1)-\mathbb{P}(Z\geq y-\beta_1)|
&\leq |\mathbb{P}(V_\lambda\geq y-\beta_1)-\mathbb{P}(Z\geq y-\beta_1)|+\mathbb{P}(W_\lambda\geq y-\beta_1,E_1^c)\\
\nonumber&\qquad+|\mathbb{P}(W_\lambda\geq y-\beta_1,E_1)-\mathbb{P}(V_\lambda\geq y-\beta_1)|\\
&\leq d_K(V_\lambda,Z)+\delta_1+|\mathbb{P}(W_\lambda\geq y-\beta_1,E_1)-\mathbb{P}(V_\lambda\geq y-\beta_1)|\,.
\end{align}
For this final term, we write
\begin{align}\label{eq:r2_3}
\nonumber|\mathbb{P}(W_\lambda\geq y-\beta_1,E_1)&-\mathbb{P}(V_\lambda\geq y-\beta_1)|\\
\nonumber&\leq|\mathbb{P}(W_\lambda\geq y-\beta_1,E_1)-\mathbb{P}(V_\lambda\geq y-\beta_1,E_1)|+\mathbb{P}(V_\lambda\geq y-\beta_1,E_1^c)\\
\nonumber&\leq|\mathbb{P}(V_\lambda\geq y-R_\lambda^\prime-\beta_1,R_\lambda^\prime\leq\beta_1)-\mathbb{P}(V_\lambda\geq y-\beta_1,R_\lambda^\prime\leq\beta_1)|+\delta_1\\
\nonumber&=\mathbb{P}(y-R_\lambda^\prime-\beta_1\leq V_\lambda<y-\beta_1,R_\lambda^\prime\leq\beta_1)+\delta_1\\
\nonumber&\leq\mathbb{P}(y-2\beta_1\leq V_\lambda<y-\beta_1)+\delta_1\\
&\leq2d_K(V_\lambda,Z)+\delta_1\,,
\end{align}
for any $\beta_1<1/2$.

Combining \eqref{eq:r2_1}--\eqref{eq:r2_3}, we thus have
\begin{equation}\label{eq:r2_final}
d_K(W_\lambda,Z)\leq3d_K(V_\lambda,Z)+\frac{8n\mu\max_jp_j}{3[4-(n+1)\max_jp_j]}\,.
\end{equation}

We now write $V_\lambda=U_\lambda+R_\lambda^{\prime\prime}$, where
\[
U_\lambda=\sum_{j=1}^rI(N_j=2)\,,\quad\text{and}\quad R_\lambda^{\prime\prime}=\sum_{j=1}^r\left(\frac{1}{p_j^\lambda\mathbb{E}[P^{-\lambda}]}-1\right)I(N_j=2)\,,
\]
and argue similarly again. For some $\beta_2>0$ to be determined later, we define the event $E_2=\{|R_\lambda^{\prime\prime}|\leq\beta_2\}$ with $\delta_2=\delta_2(\beta_2)=\mathbb{P}(E_2^c)$. We have 
\begin{equation}\label{eq:r3_1}
\delta_2=\mathbb{P}\left(\left|\sum_{j=1}^r\left(\frac{1}{p_j^\lambda\mathbb{E}[P^{-\lambda}]}-1\right)I(N_j=2)\right|>\beta_2\right)
\leq\mathbb{P}(c_\lambda U_\lambda>\beta_2)
\leq\mathbb{P}(Z>\beta_2/c_\lambda)+d_K(U_\lambda,Z)\,,
\end{equation}
where Bennett's inequality gives
\begin{align}\label{eq:r3_2}
\nonumber\mathbb{P}(Z>\beta_2/c_\lambda)
&\leq\exp\left\{-\mu\left[\frac{\beta_2}{\mu c_\lambda}\log\left(\frac{\beta_2}{\mu c_\lambda}\right)-\frac{\beta_2}{\mu c_\lambda}+1\right]\right\}\\
&=\left[\left(\frac{\beta_2}{\mu c_\lambda}\right)^{\beta_2/\mu c_\lambda}\exp\left\{1-\frac{\beta_2}{\mu c_\lambda}\right\}\right]^{-\mu}
\leq\left(\frac{e\mu c_\lambda}{\beta_2}\right)^{\beta_2/c_\lambda}\,,
\end{align}
and $\mathbb{P}(c_\lambda U_\lambda>\beta_2)$ should be interpreted as zero if $c_\lambda=0$.

Writing $d_K(V_\lambda,Z)=\sup_{y\in\mathbb{R}}|\mathbb{P}(V_\lambda\geq y-\beta_2)-\mathbb{P}(Z\geq y-\beta_2)|$, the triangle inequality gives
\begin{align}\label{eq:r3_3}
\nonumber|\mathbb{P}(V_\lambda\geq y-\beta_2)-\mathbb{P}(Z\geq y-\beta_2)|
&\leq |\mathbb{P}(U_\lambda\geq y-\beta_2)-\mathbb{P}(Z\geq y-\beta_2)|+\mathbb{P}(V_\lambda\geq y-\beta_2,E_2^c)\\
\nonumber&\qquad+|\mathbb{P}(V_\lambda\geq y-\beta_2,E_2)-\mathbb{P}(U_\lambda\geq y-\beta_2)|\\
&\leq d_K(U_\lambda,Z)+\delta_2+|\mathbb{P}(V_\lambda\geq y-\beta_2,E_2)-\mathbb{P}(U_\lambda\geq y-\beta_2)|\,.
\end{align}
For the final term here, we have
\begin{align}\label{eq:r3_4}
\nonumber|\mathbb{P}(V_\lambda\geq y-\beta_2,E_2)&-\mathbb{P}(U_\lambda\geq y-\beta_2)|\\
\nonumber&\leq
|\mathbb{P}(V_\lambda\geq y-\beta_2,E_2)-\mathbb{P}(U_\lambda\geq y-\beta_2,E_2)|+\mathbb{P}(U_\lambda\geq y-\beta_2,E_2^c)\\
\nonumber&\leq|\mathbb{P}(U_\lambda\geq y-R_\lambda^{\prime\prime}-\beta_2,|R_\lambda^{\prime\prime}|\leq\beta_2)-\mathbb{P}(U_\lambda\geq y-\beta_2,|R_\lambda^{\prime\prime}|\leq\beta_2)|+\delta_2\\
\nonumber&=\mathbb{P}(y-\beta_2-R_\lambda^{\prime\prime}\leq U_\lambda<y-\beta_2,0\leq R_\lambda^{\prime\prime}\leq\beta_2)\\
\nonumber&\qquad+\mathbb{P}(y-\beta_2\leq U_\lambda<y-\beta_2-R_\lambda^{\prime\prime},-\beta_2\leq R_\lambda^{\prime\prime}<0)+\delta_2\\
\nonumber&\leq\mathbb{P}(y-2\beta_2\leq U_\lambda<y)+\delta_2\\
&\leq2d_{K}(U_\lambda,Z)+\lfloor2\beta_2\rfloor\sup_x\mathbb{P}(Z=x)+\delta_2\,.
\end{align}

Combining \eqref{eq:r3_1}--\eqref{eq:r3_4} and choosing $\beta_2=0.49$, we have
\begin{equation}\label{eq:r3_final}
d_K(V_\lambda,Z)\leq5d_K(U_\lambda,Z)+2\left(5.55\mu c_\lambda\right)^{0.49/c_\lambda}\,.
\end{equation}

To bound $d_K(U_\lambda,Z)$ we may use Theorem 6.F of \cite{barbour92}, which gives an error bound in the Poisson approximation of the multinomial occupancy statistic $U_\lambda$. Under the assumptions that $n\geq4$ and $p_j\leq1-\sqrt{3}/2\approx0.13397$, that result gives us that
\begin{equation}\label{eq:mainterm}
d_K(U_\lambda,Z)\leq\min\left\{1,\mu^{-1}\right\}\left(6n\left[\sum_{j=1}^rp_j\pi_j\right]^2+\sum_{j=1}^r\pi_j^2+\frac{24\mu^2}{n}\right)\,.
\end{equation}

Combining \eqref{eq:r1_final}, \eqref{eq:r2_final}, \eqref{eq:r3_final} and \eqref{eq:mainterm} completes the proof of Theorem \ref{thm:main}.

\section{A generalisation of Theorem \ref{thm:main}}\label{sec:conc}

We note finally that our proof of Theorem \ref{thm:main} may be easily adapted to give a Poisson approximation bound in the following slightly more general setting: let $T=W+R$, where
\[
W=\frac{1}{g(m)}\sum_{j=1}^rh(p_j)g(N_j)\geq0\,\,\text{almost surely}\,,
\]
$\mathbb{E}[R]=0$, and $m\in\{0,1,\ldots\}$ is such that $g(0)=g(1)=\cdots=g(m-1)=0$. The power divergence statistic $\widetilde{T}_\lambda$ we consider here is the special case with $m=2$, $g=g_\lambda$, and $h$ given by $h(p_j)=(p_j\mathbb{E}[P^{-\lambda}])^{-1}$. Following the same argument as in the proof of our Theorem \ref{thm:main} but using the full generality of Theorem 6.F and Proposition A.2.5(ii) of \cite{barbour92} which we applied only for $m=2$ in the above, we have the following result.
\begin{theorem}
Let $T$ be as above and $n\geq2m$. Assume that $(1-p_j)^m\geq3/4$ and $(n+1)p_j<m+2$ for each $j=1,\ldots,r$. Let $\pi_j=\binom{n}{m}p_j^m(1-p_j)^{n-m}$ and $\mu=\pi_1+\cdots+\pi_r$, and let $Z$ have a Poisson distribution with mean $\mu$. Then
\begin{multline*}
d_K(T,Z)\leq45\min\{1,\mu^{-1}\}\left(6n\left[\sum_{j=1}^rp_j\pi_j\right]^2+\sum_{j=1}^r\pi_j^2+\frac{6m^2\mu^2}{n}\right)\\
+18(5.55\mu c)^{0.49/c}
+\frac{6(m+2)n\mu\max_jp_j}{(m+1)\left(m+2-(n+1)\max_jp_j\right)}
+8.1\text{Var}(R)\,,
\end{multline*}
where $c=\max_j|h(p_j)-1|$.
\end{theorem}   
Recalling the conditions (i), (ii), (iv) and (v) given in Section \ref{sec:intro} under which the upper bound of Theorem \ref{thm:main} goes to zero as $n\to\infty$, we note the analogous conditions here: we require $n\max_jp_j\to0$, $\text{Var}(R)\to0$, $n^m\mathbb{E}[P^{m-1}]\to\eta$ for some $0<\eta<\infty$, and $c\to0$. 

\vspace{12pt}

\noindent{\bf Acknowledgement:} The author thanks Robert Gaunt for several helpful comments, discussions and suggestions. Thanks are also due to an anonymous referee for their careful reading of the manuscript and suggestions which improved its presentation.

\end{document}